\documentclass[11pt, amsfonts]{article}

\usepackage{boxedminipage}

\def\nd{\noindent}

\def\r{{\bf R}}

\newtheorem{Teo}{Theorem}[section]

\newtheorem{lema}{Lemma}[section]
\newtheorem{defi}{Definition}[section]

\newtheorem{rmk}{Remark}[section]

\newcommand{\proof}{\normalfont{\bf Proof }}
\newcommand{\w}{W_0^{1,\Phi}}
\begin{document}

\centerline{ON VARIATIONAL MULTIVALUED  ELLIPTIC EQUATIONS ON A }
\centerline{  BOUNDED DOMAIN IN THE PRESENCE OF   CRITICAL GROWTH \footnote{Supported in part  by  CNPq/CAPES/PROCAD/UFG/UnB-Brazil}}
\bigskip

\begin{center}
{\small J. V. Goncalves~~ M. L. Carvalho }
\end{center}
\medskip

\begin{abstract}
\nd We develop arguments on the critical point theory for locally Lipschitz functionals  on Orlicz-Sobolev spaces, along with convexity and compactness techniques to investigate existence of solution of the multivalued equation
$\displaystyle  - \Delta_{\Phi} u \in \partial j(.,u) + \lambda h ~ \mbox{in}~ \Omega$, where  $\Omega \subset {\bf R}^{N}$  is a bounded smooth domain,  $\Phi : {\r} \longrightarrow [0,\infty)$  is a suitable N-function, $\Delta_{\Phi}$ is the corresponding $\Phi$-Laplacian, $\lambda > 0$ is a parameter,  $h:\Omega\rightarrow{\r}$ is integrable and $\partial j(., u)$ is the subdifferential of a  function $j$ associated with critical growth.
\end{abstract}

\centerline{\bf \scriptsize  Dedicated to Bernhard Ruf on the occasion of his $60^{th}$ birthday. }

\begin{center}
\section{Introdution}
\end{center}
\nd We deal with the multivalued equation
\begin{equation}\label{problema}
    -\Delta_\Phi u\in\partial j(.,u)+\lambda h~\mbox{in}~\Omega
\end{equation}
\nd where $\Omega \subset \r^N$ is a bounded domain with smooth boundary $\partial \Omega$, $h:\Omega\rightarrow{\r}$ is measurable, $\lambda>0$ is a parameter, $\Delta_{\Phi}$ is the
$\Phi$-Laplacian operator, that is
\[
\Delta_{\Phi} u =\mbox{div} ( \phi(|\nabla u|) \nabla u),
\]
\nd where $\phi:(0,+\infty)\rightarrow(0,+\infty)$  is continuous satisfying
$$
\begin{array}{lcl}
  (\phi_1)~~~ \mbox{(i)}~~\displaystyle  \lim_{s\rightarrow 0}s\phi (s) =0,~~ \mbox{(ii)}~~~\displaystyle  \lim_{s\rightarrow \infty} s\phi (s) = \infty,\\
  (\phi_2)~~~ s \mapsto s\phi (s)~\mbox{ is nondecreasing in}~ (0, \infty),\\
(\phi_3)~~~\mbox{there exist }~\ell,m\in(1,N)~\mbox{such that}~ \displaystyle\ell\leq \frac{t^2\phi(t)}{\Phi(t)}\leq m,~ t > 0
\end{array}
$$
\nd and $s \mapsto s \phi(s)$ is extended to ${\bf R}$ as an odd function. The functions   $\Phi,~j$  are given respectively by
\[
\Phi(t) = \int_{0}^{t}  s \phi(s) ds~~\mbox{for}~~t \in {\bf R},
\]
\begin{equation}\label{def of j}
 j(x,t) = \sigma(x)  [\Phi_*(t)-\Phi_*(a)]~\chi_{\{t > a \}}
\end{equation}
 \nd where $\sigma \in L^{\infty}(\Omega),~\sigma \geq 0,~\sigma \not\equiv 0$, $a>0$ is a number and $\Phi_*$, is the inverse of the function
$$
t\in(0,\infty)\mapsto\int_0^t\frac{\Phi^{-1}(s)}{s^{\frac{N+1}{N}}}ds
$$
\nd which extends to ${\r}$ by  $\Phi_*(t)=\Phi_*(-t)$ for  $t\leq 0$, while  $ \partial j(x,t)$ stands for the subdifferential of $j$,
\[
\partial j(x,t) = \{\mu \in  {\r}~|~  j^{o}(x,t;r) \geq \mu r,~r \in {\r}  \},
\]
\nd where $j^{o}(x,t;r)$ is the generalized directional derivative of $t \mapsto  j(x,t)$ in the direction of r,
$$
j^{o}(x,t;r) = \limsup_{y \to t,~s \to 0^{+}} \frac{j(x, y+sr) - j(x,y)}{s}.
$$
\nd Due to the nature of the differential operator $\Delta_{\Phi}$ it is natural to work in the framework of  Orlicz-Sobolev spaces.  It is known, (cf.  \cite{Fuk_1,Rao1}), that
\[
\Phi_*(t)=\int_0^t\phi_*(s)ds,
\]
\nd where $\phi_*:[0,\infty)\rightarrow[0,\infty) $  satisfies
$$
\begin{array}{lcl}
  (\phi_*)_1~~~~~~    \phi_*(0)=0,~  \phi_*(s)>0~\mbox{ for}~ s>0,~ \displaystyle \lim_{s\rightarrow\infty}\phi_*(s)=\infty,\\
  (\phi_*)_2~~~~~~     \phi_*~\mbox{is   continuous, nondecreasing},\\
(\phi_*)_3~~~~~~ \displaystyle \ell^*\leq\frac{t\phi_*(t)}{\Phi_*(t)}\leq m^*~\mbox{for}~t > 0,
\end{array}
$$
\nd  where $p^* := Np/(N-p)$ for $p \in (1,N)$. At this point we notice that
\[
j(x,t) = \int_0^t \sigma(x)~\chi_{\{\tau > a \}}~ \phi_*(\tau) d \tau,~ t \in {\bf R}.
\]
\nd  The Orlicz space associated with $\Phi$ is
$$
\displaystyle L_\Phi(\Omega):=\left\{u:\Omega\longrightarrow{\r}~\mbox{measurable}~| ~  \int_\Omega \Phi\left(\frac{u(x)}{\lambda}\right)<+\infty~ \mbox{for some}~\lambda>0\right\}
$$
 The  Orlicz-Sobolev space, (also denoted   $W^1L_\Phi(\Omega)$),  is
$$
 W^{1, \Phi}(\Omega) =\Big\{u \in L_\Phi(\Omega)~|~ \frac{\partial u}{\partial x_i} \in L_\Phi(\Omega),~ i=1,...,N \Big\}
$$
\nd and $W_0^{1,\Phi}(\Omega)$ is the closure of ${ C}_0^\infty(\Omega)$ with respect to $W^{1, \Phi}(\Omega)$.
\begin{defi}\label{solutionDef-a}  Let $h\in L_{\Phi_*}(\Omega)'$. A  vector  $u\in\w(\Omega)$ is a solution of $(\ref{problema})$ if there is an element   $\rho := \rho_u \in L_{\Phi_*}(\Omega)'$ such that
\[
\rho(x) \in\partial j(x,u(x))~\mbox{a.e.}~x\in\Omega,
\]
\[
        \int_\Omega\phi(|\nabla u|)\nabla u\nabla vdx=\int_\Omega\rho vdx+\lambda\int_\Omega hvdx,~v\in\w(\Omega).
\]
\end{defi}
\nd Our main results are,
\begin{Teo}\label{teor_principal}
Let $a > 0$ and $\ell^*>m$.  Assume that  $\phi:(0,\infty)\rightarrow(0,\infty)$ is  continuous,  satisfies $(\phi_1)-(\phi_3)$.
Let $h\in L_{\Phi_*}(\Omega)'$ be nonnegative with $h\not\equiv 0$. Then there  is   $\lambda_* >0$ such that for each $\lambda \in (0,\lambda_*)$,  equation $(\ref{problema})$ admits at least one  nonnegative solution, say $u=u_\lambda \in\w(\Omega)$.
\vskip.1cm

\nd Moreover
\begin{equation}\label{maineq ae}
-\Delta_{\Phi} u = \rho + \lambda h~~ \mbox{a.e. in}~~ \Omega.
\end{equation}
\end{Teo}

\begin{rmk}\label{example-1}
\nd {\rm If $N \geq 3$,  $\phi(t)=2$ and $\sigma \equiv 1$, then   by computing, one gets   $\Phi_*(t)= t^{\frac{2N}{N-2}}$ and $\phi_*(t) =  t^{\frac{N+2}{N-2}}$, up to constants.
\nd The subdifferential of $j(x,t)$ is shown to be
$$
\partial j(x,t)=
   \left\{
    \begin{array}{ll}
       0 ,& ~t<a \\
       \left[0,a^{\frac{N+2}{N-2}}\right], &~t=a \\
       t^{\frac{N+2}{N-2}} ,&~t>a.
    \end{array}\right.
    $$
\nd  Equation   $(\ref{problema})$ reads as
\begin{equation}\label{five}
    -\Delta u\in\partial j(.,u)+\lambda h~\mbox{in}~\Omega.
\end{equation}
\nd A nonnegative solution $u \in H_0^1(\Omega)$   of $(\ref{five})$ with $|\{x \in \Omega~|~u(x) =a \}| = 0$  is shown to satisfy
\[
-\Delta u = u^{\frac{N+2}{N-2}} \chi_{\{u > a \}} + \lambda h~~ \mbox{a.e. in}~~ \Omega.
\]
\nd Equations on bounded domains with jumping nonlinearities have been studied by many authors, see e.g. Badiale $\&$ Tarantello \cite{badiale},  Ambrosetti $\&$ Turner \cite{turner}, Chang \cite{chang}, Motreanu \& Tanaka \cite{tanaka}, Alves \& Bertone \cite{alves}   and their references.}
\end{rmk}
\nd  There is a broad literature on multivalued variational equations,  see e.g.  Halidias \&   Naniewicz \cite{halidias}, Fiacca, Matzakos \& Papageorgiou \cite{fiacca},  Alves,  Goncalves \& Santos \cite{abrantes},  Filippakis  \&  Papageorgiou \cite{filipakis}, Kyritsi  \&  Papageorgiou \cite{kyritsi},  Naniewicz \cite{naniewicz-2} and references therein.
\begin{center}
\section{Notations and Preliminary Results }
\end{center}
\nd In this section we gather  notations and results on  subdifferential calculus and Orlicz-Sobolev spaces.

\nd To begin with, following Chang \cite{chang}, Clarke \cite{Clarke1},  Motreanu \& Panagiotopoulos \cite{motreanu} and Carl,  Le \&  Motreanu \cite{CLM}, let $X$ be a reflexive real Banach space and let  $I : X \to \r$ be a locally Lipschitz continuous ($I \in Lip_{loc}(X, \r)$ for short).

\nd The generalized directional derivative of $I$ at $u \in X$ in the direction of $v\in X$ is  defined as
\[
 I^{0}(u;v) = \displaystyle \limsup_{h \to 0,~ \lambda \downarrow 0} \frac{I(u + h+ \lambda v) - I(u + h)}{\lambda}.
\]
\nd It is known that  $I^{0}(u;\cdot)$ is convex and continuous, its subdifferential at $z$ is
\[
\partial I^{0}(u;z) = \{ \mu \in X^{\prime}~ |~ \left < \mu, v - z \right > \leq I^{0}(u;v) - I^{0}(u;z)~~  v \in X \}
\]
\nd and the generalized gradient of $I$ at $u \in X$ is
\[
\partial I(u) = \{ \mu \in X^{\prime}~ |~ I^{0}(u;v) \geq  \left <\mu,v \right >,~ v\in X \}.
\]
\nd An element $u_0 \in X$ is a critical point of $I$ if  $0 \in \partial I(u_{0})$.
\vskip.1cm

\nd A main abstract result to be used in this paper is a variant for  $Lip_{loc}$ functionals,    of the Ambrosetti-Rabinowitz Mountain Pass Theorem,  to our best knowledge, developed first  via the Deformation Lemma, by  Chang \cite{chang}, see also \cite{AlvesBertoneGoncalves}  for a proof using the Ekeland Variational Principle and the Ky Fan Minimax Theorem, cf. \cite{BrezisNS}.
\vskip.1cm

\nd If $I \in Lip_{loc}(X,{\r})$ and $u \in X$ then $\partial I(u) \subset X^{\prime}$ is bounded, nonempty, convex  and weak*-closed, in the sense  that if $\xi_j\in \partial I (u_j),~ u_j \rightarrow u \mbox{ and}~  \xi_j  \stackrel {*} \rightharpoonup \xi
\mbox{ then}~ \xi\in\partial I(u)$. We set
\[
m(u) := \min_{w \in\partial I(u )} \Vert w \Vert_{X^{\prime}},~u \in X.
\]
\begin{Teo}\label{tpm}
    Let  $X$ be a  Banach space and let $I\in Lip_{loc}(X,{\r})$ with $I(0)=0$. Suppose there are numbers $\eta, r_1>0$  and  $e\in X$ such that
   $$
 \mbox{\rm (i)}~~ I(u)\geq\eta~ \mbox{if}~~\| u\|=r_1,~~ \mbox{\rm (ii)}~~\| e\|>r_1\mbox{ and}~I(e)\leq 0.
 $$
\nd Let
 \[
c =\inf_{\gamma\in\Gamma}\max_{0\leq t\leq 1}I(\gamma(t))
\]
\nd where
\[
\Gamma  =\{\gamma\in C([0,1],X)~|~\gamma(0)=0,~\gamma(1)=e\}.
\]
\nd Then $c > 0$ and there is a sequence  $(u_n)\subseteq X$ (named a $(PS)_c$-sequence) satisfying
\[
        I(u_n) \rightarrow c~~\mbox{and}~~m(u_n)\rightarrow 0.
  \]
\end{Teo}
\vskip.1cm

\nd The reader is  referred to  $\cite{A, Kuf, Rao1,gossez-Czech}$ regarding Orlicz-Sobolev spaces.  The usual norm on $L_{\Phi}(\Omega)$ is ( Luxemburg norm),
\[
\|u\|_\Phi=\inf\left\{\lambda>0~|~\int_\Omega \Phi\left(\frac{u(x)}{\lambda}\right) dx \leq 1\right\}
\]
\nd  and  the  Orlicz-Sobolev norm of $ W^{1, \Phi}(\Omega)$ is
\[
  \displaystyle \|u\|_{1,\Phi}=\|u\|_\Phi+\sum_{i=1}^N\left\|\frac{\partial u}{\partial x_i}\right\|_\Phi.
\]
\nd Recall that
$$
\widetilde{\Phi}(t) = \displaystyle \max_{s \geq 0} \{ts - \Phi(s) \},~ t \geq 0.
$$
\nd It turns out that  $\Phi$ and $\widetilde{\Phi}$  are  N-functions  satisfying  the $\Delta_2$-condition, (cf. \cite[p 22]{Rao1}).
In addition,   $L_{\Phi}(\Omega)$  and $W^{1,\Phi}(\Omega)$  are separable, reflexive,  Banach spaces. By the Poincar\'e Inequality, (see e.g.  \cite{gossez-Czech}),
\[
\int_\Omega\Phi(u)dx\leq \int_\Omega\Phi(2d|\nabla u|)dx
\]
\nd where $d=\mbox{diam}(\Omega)$, and it follows that
\[
\|u\|_\Phi\leq 2d\|\nabla u\|_\Phi~\mbox{for}~ \w(\Omega).
\]
\nd As a consequence,  $\|u\| :=\|\nabla u\|_\Phi$ defines a norm in $\w(\Omega)$, equivalent to $\|.\|_{1,\Phi}$. The imbeddings below (cf. \cite{A, Kuf, DT} ) will be  used in this paper:
\begin{equation}\label{cpt_emb}
\displaystyle \w(\Omega) \stackrel{\tiny cpt}\hookrightarrow L_\Phi(\Omega),
\end{equation}
\begin{equation}\label{Phi*Phi}
L_{\Phi^*}(\Omega) \stackrel{\tiny{cont}} \hookrightarrow L_\Phi(\Omega),
\end{equation}
\begin{equation}
W_0^{1,\Phi}(\Omega) \stackrel{\mbox{\tiny cont}}{\hookrightarrow} L_{\Phi_*}(\Omega).
\end{equation}
\nd Regarding this last case, $\Phi_*$ is the critical growth function associated  to $\Phi$, and the best constant, labeled $S$, is positive and  given by
\begin{eqnarray}\label{DT-Ineq}
\displaystyle S=\inf_{u\in\w,u \ne 0}\frac{\| u\|^\ell}{\|u\|_\Phi^\ell}.
\end{eqnarray}
\begin{rmk}\label{orlicz-sobolev}
{\rm We have  $\Delta_\Phi u\in L_{\widetilde{\Phi}_*}(\Omega~)~\mbox{for}~u \in \w$.  Indeed, set
$$
\langle -\Delta_{\Phi} u, v \rangle := \displaystyle \int_{\Omega} \phi(|\nabla u|) \nabla u \nabla v dx~~\mbox{for}~~u,v \in \w.
$$
\nd By \cite[p. 263]{Fuk_1},
$$
\displaystyle \int_\Omega\widetilde{\Phi}(\phi(|\nabla u|)|\nabla u|)dx\leq \displaystyle \int_\Omega \Phi(2|\nabla u|)dx<\infty,
$$
\nd which gives $\phi(|\nabla u|)|\nabla u|\in L_{\widetilde{\Phi}}(\Omega)$. By the H\" older inequality,
$$
    \displaystyle \int_{\Omega} |\phi(|\nabla u|) \nabla u \nabla v| dx \leq 2 \|\phi(|\nabla u|)|\nabla u| \|_{\widetilde{\Phi}}\| v\|.
$$
\nd As a consequence of the inequality above and $(\ref{Phi*Phi})$,   $\Delta_\Phi u\in L_{{\Phi}}(\Omega~)^{\prime} = L_{\widetilde{\Phi}}(\Omega~) \subseteq L_{\widetilde{\Phi}_*}(\Omega~)$.
}
\end{rmk}
\nd The energy functional associated with  $(\ref{problema})$ is $I := I_{\lambda,a}$ defined by
\[
I(u)=\int_\Omega\Phi(|\nabla u|)dx-\int_\Omega j(x,u)dx-\lambda\int_\Omega hudx,~u\in\w.
\]
\nd  Set
\[
Q_{\lambda}(u) =\int_\Omega\Phi(|\nabla u|)dx-\lambda\int_\Omega hudx~~\mbox{and}~~ J(u):=\int_\Omega j(x,u)dx.
\]
\nd It is known that
\[
Q_{\lambda}\in {C}^1(\w, {\r}),~~
\langle Q_{\lambda}^{\prime}(u),v \rangle = \displaystyle \int_\Omega\phi(|\nabla u|)\nabla u\nabla vdx-\lambda\int_\Omega hvdx,
\]
\noindent and, (cf. lemma \ref{aubin_clarke_Le}),
\begin{equation}\label{aubin-clark}
J\in Lip_{loc}(\w,{\r})~\mbox{and}~\partial J(u)\subseteq \Big\{\rho \in L_{\widetilde{\Phi_*}}(\Omega)~|~\rho(x)\in\partial j(x,u(x))~\mbox{a.e.}~ x\in\Omega \Big\}.
\end{equation}
\nd Thus,
$$
I \in Lip_{loc}(\w(\Omega),{\r})~\mbox{and}~\partial I(u)=Q_{\lambda}'(u)- \partial J(u).
$$
\nd Moreover,  $u$ is a critical point of $I$ if $0 \in \partial I(u)$ that is,  there is some $\rho\in \partial J(u) $~such that
\[
\langle Q_{\lambda}^{\prime}(u),v \rangle - \langle \rho, v \rangle = 0~\mbox{for}~v \in \w.
\]
\begin{center}
\section{The Mountain Pass Geometry of $I$}
\end{center}

\nd The proof of theorem \ref{teor_principal} uses  theorem \ref{tpm}. Items (i)-(ii) in theorem \ref{tpm} are known as the mountain pass geometry for $I$. In this regard we will present a proof of the result below based on  \cite{alves}.
\vskip.2cm

\begin{lema}\label{gpm}
   Let $h\in L_{\widetilde{\Phi}_*}(\Omega)$ be nonnegative, with $h \ne 0$, and assume that  $\ell^*>m$. Then there exist  $\lambda_0, \eta, r_1>0$ and  $e\in\w$ such that for each $\lambda\in(0,\lambda_0)$ and $a > 0$,
    \begin{description}
      \item{\rm (i)}~~~ $I(u)\geq\eta>0~\mbox{if}~\| u\|=r_1,$
      \item{\rm (ii)}~~~ $\|e\|>r_1~\mbox{and}~I(e) \leq 0$.
    \end{description}
    \end{lema}
\nd \proof ~ At first we show  {\rm (i)}. Indeed, using  lemmas   \ref{lema_naru}, (cf. Appendix), and the H\"older Inequality we have
\begin{equation}\label{for I}
 I(u)  \geq  \displaystyle \min\{\| u\|^\ell,\| u\|^m\}-\int_\Omega j(x,u)dx-2\lambda\|h\|_{\widetilde{\Phi}_*}\|u\|_{\Phi_*}.
\end{equation}
\nd Using lemma \ref{lema_naru_1} also in the  Appendix, we get

\begin{equation}\label{j-ESTIMATE}
\begin{array}{lcl}
   \displaystyle \int_\Omega j(x,u)dx + 2\lambda\|h\|_{\widetilde{\Phi}_*}\|u\|_{\Phi_*} \leq\\ \\
    |\sigma|_\infty\max\{\|u\|_{\Phi_*}^{\ell^*},\|u\|_{\Phi_*}^{m^*}\}+2\lambda\|h\|_{\widetilde{\Phi}_*}\|u\|_{\Phi_*} \leq\\ \\
    |\sigma|_\infty\max\left\{\frac{1}{S^{\frac{\ell^*}{\ell}}}\| u\|^{\ell^*},\frac{1}{S^{\frac{m^*}{\ell}}}\| u\|^{m^*}\right\}+\frac{2\lambda}{S^{\frac{1}{\ell}}}\|h\|_{\widetilde{\Phi}_*}\| u\|.\\ \\
\end{array}
\end{equation}
\nd Joining estimates (\ref{for I}) and (\ref{j-ESTIMATE}) we have
$$
I(u)  \geq  \displaystyle \min\{\| u\|^\ell,\| u\|^m\} - |\sigma|_\infty\max\left\{\frac{1}{S^{\frac{\ell^*}{\ell}}}\| u\|^{\ell^*},\frac{1}{S^{\frac{m^*}{\ell}}}\| u\|^{m^*}\right\}-\frac{2\lambda}{S^{\frac{1}{\ell}}}\|h\|_{\widetilde{\Phi}_*}\| u\|.
$$
\nd Taking   $\| u\| \leq 1$ it follows by the inequality just above that
$$
I(u)\geq \|u\|^m\left(1-\beta\| u\|^{\ell^*-m}-\alpha\lambda\|u\|^{1-m}\right),
$$
\nd  where $\alpha:= {2}/{S^{\frac{1}{\ell}}}\|h\|_{\widetilde{\Phi}_*},˜~\beta:={|\sigma|_\infty}/{S^{\frac{\ell^*}{\ell}}}$.  Set
$P(s):=1-\beta s^{\ell^*-m},~s > 0$. Since  $\ell^*>m$ one gets
$$
\begin{array}{lll}
P(s) \geq {1}/{2}~ \mbox{if}~ s \leq s_0:=\left(\frac{1}{2\beta}\right)^{{1}/{\ell^{*}-m}},\\
\\
P(s_0)-\lambda\alpha s_0^{1-m}\geq \frac{1}{4}~\mbox{whenever}~  \lambda \leq \lambda_0:=\frac{1}{4\alpha}s_0^{m-1}.
\end{array}
$$
\nd Choosing  $r_1 := \min\{1,s_0\}$ it follows that
$$I(u)\geq\frac{r_1^m}{4}>0~\mbox{for}~u\in\w~\mbox{with}~\| u\|=r_1.$$
\nd This shows $\mbox{(i)}$. In order to show (ii), pick $\varphi\in\mathbf{C}^\infty_0(\Omega)$ with $\varphi \geq 0$ such that
$$
\mbox{meas}\{x\in\Omega~|~\varphi(x)\geq a\}>0~~\mbox{and}~~\|\varphi\|\geq 1.
$$
\nd Taking $t>1$  we get
$$
\begin{array}{lll}
   I(t\varphi) &\leq&   t^m\| \varphi\|^m-  \int_{\{\varphi\geq a\}}\sigma(x)\Phi_*(t\varphi)dx+|\sigma|_\infty\Phi_*(a)|\Omega| \\ \\
   & \leq &  t^m\|\varphi\|^m-t^{\ell^*}\int_{\{\varphi\geq a\}}\sigma(x)\Phi_*(\varphi)+|\sigma|_\infty\Phi_*(a)|\Omega|.
\end{array}
$$
\nd As a consequence,
\[
I(t\varphi)\stackrel{t\rightarrow\infty}\longrightarrow-\infty.
\]
\nd Setting  $e:=t_1\varphi$ with $t_1>1$ large enough we have $I(e)<0$, showing (ii).   $\hfill{\rule{2mm}{2mm}}$
\begin{center}
\section{Boundedness of the Palais-Smale Sequence}
\end{center}
\nd The  result  below is a  special case of theorem 1.1 in  Le, Motreanu and Motreanu \cite{Le2} which in turn  is a  variant for Orlicz-Sobolev spaces, of the Aubin-Clarke Theorem  (cf. \cite[theorem 2.7.5]{Clarke1}). The result itself as well as its proof will be used several times in this paper.
\vskip.2cm

\begin{lema}\label{aubin_clarke_Le}
    Let $j$ be as in $(\ref{def of j})$. Then
\[
\partial j(x,t)=
\left\{
    \begin{array}{ll}
       0,~~t<a,  \\
       \left[0,\sigma(x)\phi_*(a)\right],~t=a,  \\
       \sigma(x)\phi_*(t),~t>a,
\end{array}\right.
\]
\nd  and the functional
    $$J(u):=\int_\Omega j(x,u(x))dx,~u\in L_{\Phi_*}(\Omega)$$
\nd satisfies
$$
\begin{array}{lcl}
J  \in Lip_{loc}( L_{\Phi_*}(\Omega), {\bf R})~ \\
\nd \mbox{and}\\
\partial J(u)\subseteq\{\rho \in L_{\widetilde{\Phi_*}}(\Omega)~|~\rho(x)\in\partial j(x,u(x))\mbox{ a.e. }x\in\Omega\}.
\end{array}
$$
\end{lema}

\nd By  lemmas \ref{gpm}, \ref{aubin_clarke_Le} and theorem \ref{tpm} there is a sequence $(u_n)\subseteq\w$ such that
\begin{equation}\label{cond_ps}
        I(u_n)\stackrel{n}\longrightarrow c~~\mbox{and}~~m(u_n) \equiv \min_{w\in\partial I(u_n)}\|w\|_{W^{-1,\widetilde{\Phi}}}\stackrel{n}\longrightarrow 0.
\end{equation}
\nd Actually, there is $w_n \in \partial I(u_n)$ such that
$$
\|w_n\|_{W^{-1,\widetilde{\Phi}}} = \min_{w\in\partial I(u_n)}\|w\|_{W^{-1,\widetilde{\Phi}}}
$$
\nd and so  there is $\rho_n \in \partial J(u_n)$ such that $w_n = Q_{\lambda}^{\prime}(u_n) - \rho_n$. Hence, $\langle w_n, v \rangle \rightarrow 0,~ v \in \w(\Omega)$

\nd so that

\begin{equation}\label{e1}
   \displaystyle   \int_{\Omega}\phi(|\nabla u_n|)\nabla u_n\nabla v dx =  \lambda\displaystyle\int_{\Omega}h v dx+\int_{\Omega}\rho_n v dx + o_n(1).
    \end{equation}

\nd The result below is inspired on lemma 1.20 of Willem \cite{Willem}.
\vskip.2cm

\begin{lema}\label{seq_lim}
\nd The $(PS)_{c}$ - sequence  $(u_n)\subseteq\w$  is bounded. In particular,  there is some  $u^1\in\w$ such that
    $$\displaystyle u_n\rightharpoonup u~\mbox{in }\w.$$
\end{lema}
\proof {\bf of Lemma \ref{aubin_clarke_Le}} By the very definition of $j$, $j(x,.)$ is differentiable at each  $t\neq a$ and
\[
\partial j(x,t)=j'(x,t)= \sigma(x)\phi_*(t) \chi_{\{ t > a \}}.
\]
\nd On the other hand, if $t=a$, then (cf. \cite{CLM}),
\[
\partial j(x, a)=\left[\lim_{t\rightarrow a^-}~\chi_{\{t > a \}}\sigma(x)\phi_*(t),\lim_{t\rightarrow a^+}~\chi_{\{t > a \}} ~\sigma(x)\phi_*(t)\right]=[0,\sigma(x)\phi_*(a)].
\]
\nd In particular,  for each  $\rho = \rho(x) \in\partial j(x,t)$ with $t\geq0$ we have,
\[
    0\leq t\rho \leq \sigma(x)t\phi_*(t)\leq m^*\sigma(x)\Phi_*(t)~~\mbox{a.e.}~~x\in\Omega.
\]
\nd Actually, if  $t>a$, then
\begin{equation}\label{Jul_03-6}
    \ell^*\sigma(x)\Phi_*(t)\leq t\rho\leq m^*\sigma(x)\Phi_*(t)~~\mbox{a.e.}~~x\in\Omega.
\end{equation}
\nd {\rm Notice that  if  $\rho := \rho(x) \in\partial j(x,t)$ then
    $$0\leq\rho\leq\sigma(x)\phi_*(t)\leq|\sigma|_\infty\phi_*(t).$$
 \nd Moreover, using the fact that  $\widetilde{\Phi}_*(\phi_*(t))\leq\Phi_*(2t)$,  (cf. \cite[p. 263]{Fuk_1}), we infer that
    $$
0\leq\rho\leq|\sigma|_\infty\widetilde{\Phi}_*^{-1}\circ\Phi_*(2t),
$$
\nd which is condition (1.6) in  theorem 1.1 of \cite{CLM}. This proves lemma 5.1.}$\hfill{\rule{2mm}{2mm}}$
\vskip.1cm

\nd \proof {\bf of Lemma \ref{seq_lim}.} By (\ref{cond_ps})   we have
\[
    |\langle w_n,u_n\rangle|\leq m(u_n) \|u_n\|\leq \ell^*\|u_n\|~ \mbox{for}~n~ \mbox{large enough}.
\]
\nd Set
$$
S_I:=\sup_{n}I(u_n)<\infty.
$$
\nd Estimating using  the inequality above,  (\ref{Jul_03-6}), the H\"older Inequality and lemma \ref{lema_naru}  we have
$$
\begin{array}{lll}
  S_I+\|u_n\| & \geq & \displaystyle I(u_n)-\frac{1}{\ell^*}\langle w_n,u_n\rangle \\ \\
  & \geq & \displaystyle\left(1-\frac{m}{\ell^*}\right)\int_\Omega\Phi(|\nabla u_n|)dx -\lambda\left(1-\frac{1}{\ell^*}\right)\int_\Omega h u_ndx+\frac{1}{\ell^*}\int_{\{u_n= a\}}\rho_nadx\\ \\
  & & +\displaystyle\int_{\{u_n> a\}}\left[\frac{1}{\ell^*}\rho_nu_n-j(x,u_n)\right]dx \\ \\
  & \geq & \displaystyle\left(1-\frac{m}{\ell^*}\right)\int_\Omega\Phi(|\nabla u_n|)dx+\int_{\{u_n>a\}}[\sigma(x) \Phi_*(u_n)-j(x,u_n)]dx\\ \\
 & &\displaystyle  -\lambda\left(1-\frac{1}{\ell^*}\right)\int_\Omega h u_ndx \\ \\
  & \geq & \displaystyle\left(1-\frac{m}{\ell^*}\right)\int_\Omega\Phi(|\nabla u_n|)dx -2\lambda\left(1-\frac{1}{\ell^*}\right)\|h\|_{\widetilde{\Phi}}\|u_n\|_{\Phi^*}\\ \\
  & \geq & \displaystyle\left(1-\frac{m}{\ell^*}\right)\min\{\|u_n\|^\ell,\| u_n\|^m\}-\frac{2\lambda}{S^{\frac{1}{\ell}}}\left(1-\frac{1}{\ell^*}\right)\|h\|_{\widetilde{\Phi}}\|u_n\|,
\end{array}
$$
\nd showing that  $(\|u_n\|)$ is bounded. $\hfill{\rule{2mm}{2mm}}$
\begin{center}
\section{ On the Convergence of the Palais-Smale Sequence}
\end{center}
\nd  The result below is crucial, will be proved in detail in this paper, and actually, was motivated by lemma 4.4 by Fukagai,  Ito \& Narukawa \cite{Fuk_1} .
\vskip.1cm

\begin{lema}\label{conv_compactos}
    Let  $(u_n) \subset \w(\Omega)$ be the sequence in $(\ref{cond_ps})$. Extend each $u_n$ to ${\bf R}^{N}$
   by setting $u_n = 0~\mbox{on}~ {\bf R}^{N} \backslash \Omega$. Then there are $x_1, \cdots, x_r \in {\bf R}^{N}$ such that
 \begin{equation}\label{conv-K}
    u_n\stackrel{L_{\Phi_*}(K)}\longrightarrow u
    \end{equation}
    for each compact set $K\subset \r^N \backslash \{x_1, \cdots , x_r \}$.
\end{lema}
\nd At first we gather some notatios and remarks, (cf.  Willem \cite{Willem}). Given $v\in {C}^\infty_0(\Omega)$ we extend it to  $\r^N$ by setting $v(x)=0$ if $x\in \r^N \backslash \Omega$ and denote the extension by $v$. Then $v\in {C}_0^\infty(\r^N)$ and $\mbox{supp}(v)\subseteq \Omega $. In addition,
$$
\| v\|_{W^{1,\Phi}(\r^N)}=\| v\|_{W^{1,\Phi}(\Omega)}
$$
\nd and
$$
\w(\Omega)=\overline{\{v\in {C}_0^\infty(\r^N)~|~\mbox{supp}(v)\subseteq\Omega\}}^{W^{1,\Phi}(\r^N)}.
$$
\nd Thus, if $v\in\w(\Omega)$ then $v \in {W^{1,\Phi}(\r^N)}$. Similar notations for functions in $L_{\widetilde{\Phi}_*}(\Omega)$.
\vskip.1cm

\nd Consider the normed space
$$
{C}_0 =\overline{\{u\in {C}(\Omega)~|~\mbox{supp}(u) \buildrel \mbox{\scriptsize {cpt}} \over \subseteq \r^N\}}^{|\cdot|_\infty},
$$
\nd where $\displaystyle |u|_\infty=\sup_{x\in \r^N}|u(x)|$ and denote by $\mathcal{M}$ the space of finite measures on $\r^N$  with the norm
$$
\|\mu\|_{\mathcal{M}}=\sup\left\{\int u d\mu~|~u\in {C}_0,~ |u|_\infty=1\right\}.
$$
\begin{rmk}
We recall below some  notations and results:
\begin{description}
\item{\rm(i)} $\mathcal{M} = {C}_0^*$  and $\langle \mu,u \rangle = \int u d\mu$,

\item{\rm(ii)}  $\mu_n \buildrel \mathcal{M} \over \rightharpoonup \mu$ means that $\int u d\mu_n\stackrel{n\rightarrow\infty}{\longrightarrow} \int u d\mu,~u\in {C}_0$,

\item{\rm(iii)} if $(\mu_n)\subseteq \mathcal{M}$ is bounded then $\mu_n \buildrel \mathcal{M} \over \rightharpoonup \mu$, up to subsequence.
\end{description}
\end{rmk}
\nd By lemma \ref{seq_lim} the  $(PS)_c$-sequence $(u_n)\subseteq \w(\Omega)$ is bounded.  Consider $\mu_n,\nu_n: {C}_0\rightarrow \r$,
$$
\langle\mu_n, v\rangle=\int_{\r^N}\Phi(|\nabla u_n|)vdx~~\mbox{and}~~\langle\nu_n, v\rangle=\int_{\r^N}\Phi_*(|u_n|)vdx,~v\in {C}_0.
$$
\nd Then there is a constant $C > 0$ such that
$$
|\langle\mu_n,v\rangle| \leq C |v|_\infty~\mbox{ and}~  |\langle\nu_n,v\rangle|\leq C|v|_\infty
$$
\nd that is $(\mu_n),(\nu_n)\subseteq \mathcal{M}$ are bounded. It follows that
\begin{equation}\label{conv_med}
\Phi(|\nabla u_n|)\rightharpoonup \mu,~~\Phi_*(|u_n|)\rightharpoonup \nu~~\mbox{in}~~\mathcal{M}.
\end{equation}
\nd We shalll need the following variant for Orlicz-Sobolev spaces of the concentration-compactness principle cf. Lions \cite{lions1}, Fukagai, Ito \& Narukawa  \cite{Fuk_1}.
\begin{lema}\label{conc_comp}
\nd  There exist a denumerable set $J$, a family $\{x_j\}_{j\in J}\subseteq {\r}^N$  with $x_i \ne x_j$ and families of nonnegative numbers $\{\nu_j\}_{j\in J}$ and $\{\mu_j \}_{j\in J}$ such that
$$
\nu=\Phi_*(u^1)+\sum_{j\in J}\nu_j\delta_{x_j}~\mbox{and}~\mu\geq \Phi(|\nabla u^1|)+\sum_{j\in J}\mu_j\delta_{x_j},
$$
\nd where  $\delta_{x_j}$ is the Dirac measure with mass at $x_j$. In addition,
$$
\nu_j\leq\max\left\{S^{-\frac{\ell^*}{\ell}}\mu_j^{\frac{\ell^*}{\ell}},S^{-\frac{m^*}{\ell}}\mu_j^{\frac{m^*}{\ell}},
      S^{-\frac{\ell^*}{\ell}}\mu_j^{\frac{\ell^*}{m}},S^{-\frac{m^*}{\ell}}\mu_j^{\frac{m^*}{m}}\right\},~j\in J.
$$
\end{lema}
\begin{lema}\label{J_finito}
   The  set $\widetilde{J} = \{j\in J~|~\nu_j>0\}$ is finite.
\end{lema}
\nd \proof  We claim that  $\{x_j\}_{j\in \widetilde{J}}\subseteq \overline{\Omega}$. Indeed, if on the contrary,  $x_j\in\overline{\Omega}^c$ for some $j\in \widetilde{J}$, there is $\epsilon > 0$ such that ${\overline{B}_\epsilon}(x_j)\subseteq\overline{\Omega}^c$. Choose $\varphi_\epsilon\in {C}_0^\infty(\r^N)$ such that
$$
supp(\varphi_\epsilon)\subseteq B_\epsilon(x_j),~~~
\varphi_\epsilon\stackrel{\epsilon\rightarrow 0}\longrightarrow \chi_{\{x_j\}}~ \mbox{a.e.}~ {\bf R}^{N}.
$$
\nd Now, we extend  $u_n$ to $\r^N$ by setting $u_n(x)=0$ for $x\in\r^N-\Omega$. Take $\epsilon > 0$. Using  $(\ref{conv_med})$, we have
$$
0=\int_{\r^n}\Phi(|\nabla u_n|)\varphi_\epsilon dx\stackrel{n} \longrightarrow \int_{\r^N}\varphi_\epsilon d\mu,
$$
\nd and passing to the limit as $\epsilon \to 0$ we get,
$$
0=\int_{\r^N}\varphi_\epsilon d\mu=\int_{B_\epsilon(x_j)}\varphi_\epsilon d\mu  \rightarrow\int_{\{x_j\}}d\mu=\mu_j.
$$
\nd Hence, $\mu_j = 0$ and by lemma \ref{conc_comp} we infer that $\nu_j=0$, impossible because $j\in\widetilde{J}$, showing the claim.
\vskip.1cm

\nd We claim that
\[
(\phi(|\nabla u_n|) |\nabla u_n|)~ \mbox{is bounded in}~ L_{\widetilde{\Phi}}(\Omega)
\]
\nd Indeed, take  $\psi\in C_0^\infty$ such that $0 \leq \psi \leq 1$, $\psi(x) = 1~ \mbox{if}~ |x| \leq 1$ and $\psi(x) = 0~\mbox{if}~|x| \geq 2$. Pick $x_j$ with $j\in\widetilde{J}$, $\epsilon > 0$ and set
$$
\psi_\epsilon(x):=\psi\left(\frac{x-x_j}{\epsilon}\right),~~x\in{\r}^N.
$$
\nd Notice that  $(\psi_\epsilon u_n)\subseteq\w(\Omega)$ is bounded. At this point we recall that
\begin{equation}\label{MM}
w_n = Q_{\lambda}^{\prime}(u_n) - \rho_n~~ \mbox{for some}~~\rho_n \in \partial J(u_n).
\end{equation}
\nd Since $m(u_n) \rightarrow 0$ we infer from (\ref{MM}) that
\begin{equation}
    \begin{array}{lll}\label{e2}
      \displaystyle\int_{\Omega}\phi(|\nabla u_n|)\nabla u_n\nabla (\psi_\epsilon u_n) & = & \lambda\displaystyle\int_{\Omega}hu_n\psi_\epsilon dx+\int_{\Omega}\rho_nu_n\psi_\epsilon dx+o_n(1).
    \end{array}
\end{equation}
\nd Moreover, by lemma \ref{aubin_clarke_Le},  $\rho_n\in L_{\widetilde{\Phi}_*}(\Omega)$ and $\rho_n(x) \in\partial j(x,u_n(x))$ for $x\in\Omega$. By (\ref{e2}) and  lemma \ref{aubin_clarke_Le},
\begin{equation}\label{d1}
    \begin{array}{lll}
        \displaystyle\int_{\Omega}\phi(|\nabla u_n|)\nabla u_n\nabla (\psi_\epsilon u_n) & = & \displaystyle\left(\int_{\{u_n< a\}}+\int_{\{u_n\geq a\}}\right)\rho_nu_n\psi_\epsilon dx+\lambda\int_{\Omega}hu_n\psi_\epsilon dx+o_n(1)\\ \\
        & = & \displaystyle\int_{\{u_n\geq a\}}\rho_nu_n\psi_\epsilon dx+\lambda\int_{\Omega}hu_n\psi_\epsilon dx+o_n(1)\\ \\
        & \leq &  \displaystyle m^*\int_{\{u_n\geq a\}} \sigma(x)\Phi_*(u_n)\psi_\epsilon dx+\lambda\int_{\Omega}hu_n\psi_\epsilon dx+o_n(1)\\ \\
        & \leq & \displaystyle m^*|\sigma|_\infty\int_\Omega\Phi_*(u_n)\psi_\epsilon dx+\lambda\int_{\Omega}hu_n\psi_\epsilon dx+o_n(1)
    \end{array}
\end{equation}
\nd On the other hand, using the fact that $t^2\phi(t)\geq\Phi(t)$ we have,
\begin{equation}\label{d2}
    \begin{array}{lll}
        \displaystyle\int_{\Omega}\phi(|\nabla u_n|)\nabla u_n\nabla (\psi_\epsilon u_n) & = & \displaystyle\int_{\Omega}u_n\phi(|\nabla u_n|)\nabla u_n\nabla \psi_\epsilon dx+\int_{\Omega}\psi_\epsilon\phi(|\nabla u_n|)|\nabla u_n|^2dx \\ \\
        & \geq & \displaystyle\int_{\Omega}u_n\phi(|\nabla u_n|)\nabla u_n\nabla \psi_\epsilon dx+\int_{\Omega}\psi_\epsilon\Phi(|\nabla u_n|)dx
    \end{array}
\end{equation}
\nd Using (\ref{d1}), (\ref{d2}) and the inequality  $\widetilde{\Phi}(t\phi(t))\leq\Phi(2t)$ it follows that
$(\phi(|\nabla u_n|)|\nabla u_n|)$ is bounded in $L_{\widetilde{\Phi}}(\Omega)$, showing the claim.
\vskip.1cm

\nd As a consequence
$(\phi(|\nabla u_n|){\partial u_n}/{\partial x_i})$ is also bounded in $L_{\widetilde{\Phi}}(\Omega)$ and so
\begin{equation}\label{e4}
    \displaystyle \phi(|\nabla u_n|)\frac{\partial u_n}{\partial x_i}\rightharpoonup w_i~\mbox{in}~L_{\widetilde{\Phi}}(\Omega),~i=1,...,N.
\end{equation}
\nd  Setting  $w =(w_1,...,w_N)$, we claim that
\begin{equation}\label{e6}
    \int_\Omega (u_n\phi(|\nabla u_n|)\nabla u_n\nabla\psi_\epsilon-u~ w.\nabla \psi_\epsilon)dx=o_n(1).
\end{equation}
\nd Indeed, in a first step applying  an easy estimate  and in a second step  using the the H\"older inequality and  applying  (\ref{e4}) with test function  $\frac{\partial \psi_\epsilon}{\partial x_i}u$, we have
$$
\begin{array}{cl}
  \displaystyle  \Big |\int_\Omega  \phi(|\nabla u_n|)\frac{\partial u_n}{\partial x_i}\frac{\partial \psi_\epsilon}{\partial x_i}u_n-  w_i\frac{\partial \psi_\epsilon}{\partial x_i}u  dx \Big |~  \leq\\
\\
\displaystyle\int_\Omega \Big | \phi(|\nabla u_n|)\frac{\partial u_n}{\partial x_i}\frac{\partial \psi_\epsilon}{\partial x_i}\big(u_n- u \big) \Big | dx+
        \Big |  \int_\Omega \phi(|\nabla u_n|)\frac{\partial u_n}{\partial x_i}\frac{\partial \psi_\epsilon}{\partial x_i}u - w_i\frac{\partial \psi_\epsilon}{\partial x_i}u  dx \Big| ~ \leq \\
\\
\displaystyle 2\left\|\phi(|\nabla u_n|)\frac{\partial u_n}{\partial x_i}\frac{\partial \psi_\epsilon}{\partial x_i}\right\|_{\widetilde{\Phi}}\left\|u_n-u\right\|_\Phi+ o_n(1)
\end{array}
$$
\nd  Since by  $(\ref{cpt_emb})$, $\|u_n-u\|_\Phi \rightarrow 0$, we infer that
\[
    \int_\Omega \phi(|\nabla u_n|)\frac{\partial u_n}{\partial x_i}\frac{\partial \psi_\epsilon}{\partial x_i}u_ndx\stackrel{n}\longrightarrow\int_\Omega w_i\frac{\partial \psi_\epsilon}{\partial x_i}u dx,~~i=1,...,N,
\]
\nd which leads to (\ref{e6}),  showing the claim. Replacing  (\ref{e6}) in (\ref{d2}) we get
\begin{equation}\label{d3}
    \begin{array}{lll}
	\displaystyle\int_{\Omega}\psi_\epsilon\Phi(|\nabla u_n|)dx
        +\int_\Omega u w.\nabla \psi_\epsilon dx  & \leq &
        \displaystyle\int_{\Omega}\phi(|\nabla u_n|)\nabla u_n\nabla (\psi_\epsilon u_n)+o_n(1).
    \end{array}
\end{equation}
\nd It follows from (\ref{d1}) and (\ref{d3}) that
\[
    \displaystyle\int_{\Omega}\psi_\epsilon\Phi(|\nabla u_n|)dx+\int_\Omega u w.\nabla \psi_\epsilon dx\leq m^*|\sigma|_\infty\int_\Omega \Phi_*(u_n)\psi_\epsilon dx+\lambda\int_{\Omega}hu_n\psi_\epsilon dx+o_n(1).
\]
\nd Passing to the limit in the inequality just above in $n$, recalling that
$$
\int_{\Omega}\Phi(|\nabla u_n|)\psi_\epsilon dx\stackrel{n\rightarrow\infty}\longrightarrow \int_{\Omega}\psi_\epsilon d\mu,~~~ \int_{\Omega}\Phi_*( u_n)\psi_\epsilon dx\stackrel{n\rightarrow\infty}\longrightarrow \int_{\Omega}\psi_\epsilon d\nu
$$
\nd and
$$
\int_{\Omega}u_nh\psi_\epsilon dx\stackrel{n\rightarrow\infty}\longrightarrow\int_{\Omega}u h\psi_\epsilon dx.
$$
\nd we get to
\begin{equation}\label{d5}
    \displaystyle\int_{\Omega}\psi_\epsilon d\mu+\int_\Omega u w.\nabla \psi_\epsilon dx\leq m^*|\sigma|_\infty\int_{\Omega}\psi_\epsilon d\nu+\lambda\int_{\Omega}hu \psi_\epsilon dx.
\end{equation}
\nd We claim that $(\rho_n)$ is bounded in  $L_{\widetilde{\Phi}_*}(\Omega)$. Indeed, using lemma \ref{aubin_clarke_Le}
we get
$$
\begin{array}{lll}
\displaystyle\int_\Omega\widetilde{\Phi}_*(\rho_n)dx & \leq & \displaystyle\int_\Omega\widetilde{\Phi}_*(\sigma(x)\phi_*(u_n))dx\leq \int_\Omega\widetilde{\Phi}_*(|\sigma|_\infty\phi_*(u_n))dx \\ \\
& \leq & \displaystyle C_{|\sigma|_\infty} \int_\Omega\widetilde{\Phi}_*(\phi_*(u_n))dx\leq C_{|\sigma|_\infty}\int_\Omega\Phi_*(2u_n)dx\leq C,
\end{array}
$$
\nd showing the claim. Thus there is    $\rho\in L_{\widetilde{\Phi}_*}(\Omega)$ such that
$$
\rho_n\rightharpoonup\rho~\mbox{in}~L_{\widetilde{\Phi}_*}(\Omega).
$$
\nd Let $v\in\w(\Omega)$. Passing to the limit in the expression
$$
\langle w_n,v\rangle=\int_\Omega\left(\phi(|\nabla u_n|)\nabla u_n\nabla v-\lambda hv -\rho_n v\right)dx,
$$
\nd and using (\ref{e4}) we get to
\begin{equation}\label{d6}
\int_\Omega\left(w.\nabla v-\lambda hv -\rho v\right)dx=0.
\end{equation}
\nd Setting  $v=u \psi_\epsilon$ in (\ref{d6}) we have
\[
    \int_\Omega u  w.\nabla \psi_\epsilon dx=\int_\Omega (\lambda hu +\rho u - w.\nabla u)\psi_\epsilon dx.
\]
\nd But
$$
|(\lambda hu +\rho u-w.\nabla u)\psi_\epsilon|\leq |\lambda hu|+|\rho u|+|w.\nabla u|\in L^1(\Omega)
$$
\nd and
$$
(\lambda hu +\rho u-w.\nabla u)\psi_\epsilon\stackrel{\epsilon\rightarrow 0}\longrightarrow 0~\mbox{a.e. in}~\Omega
$$
\nd By means of Lebesgue's theorem,
$$
\int_\Omega  u w.\nabla \psi_\epsilon dx\stackrel{\epsilon\rightarrow 0}\longrightarrow 0~\mbox{and}~ \int_\Omega hu \psi_\epsilon dx\stackrel{\epsilon\rightarrow 0}\longrightarrow 0.
$$
\nd Noticing that
$$
\psi_\epsilon \stackrel{\epsilon\rightarrow 0}\longrightarrow \chi_{\{x_j\}}~\mbox{a.e. in}~ \r^N~\mbox{and}~ \psi_\epsilon(x)\leq \chi_{B_1(x_j)}(x)~\mbox{for}~x \in {\bf R}^{N},~ \epsilon>0~ \mbox{small}
$$
\nd we get to
$$
\int_{\r^N} \psi_\epsilon d \mu\stackrel{\epsilon\rightarrow 0}\longrightarrow \int_{\{x_j\}}d\mu=\mu(\{x_j\})=\mu_j~\mbox{and}~ \int_{\r^N} \psi_\epsilon d \nu\stackrel{\epsilon\rightarrow 0}\longrightarrow \int_{\{x_j\}}d\nu=\nu(\{x_j\})=\nu_j.
$$
\nd Passing to the limit in (\ref{d5}) we get to
\begin{equation}\label{d8}
    \mu_j\leq m^*|\sigma|_\infty\nu_j,~~j\in \widetilde{J}.
\end{equation}
\nd By lemma \ref{conc_comp}, $\mu_j\leq c_1\mu_j^\alpha$, where $1<\alpha\leq\min \big\{{\ell^*}/{\ell},{m^*}/{\ell},{\ell^*}/{m},{m^*}/{m} \big \}.$
\vskip.1cm

\nd Thus $\mu_j \geq c_2$ for some positive constant $c_2$. In addition by (\ref{d8}), $\nu_j\geq c_3,~\mbox{for}~j\in\widetilde{J}$ and for some positive constant $c_3$. At this point, we infer that if $\#(\widetilde{J}) = \infty$, then
$$
\sum_{j\in\widetilde{J}}\nu_j\geq\sum_{j\in\widetilde{J}}c_3=\infty,
$$
\nd which is impossible because $\nu$ is a finite measure and
$$
\nu=\Phi_*(u)+\sum_{j\in \widetilde{J}}\nu_j\delta_{x_j}.
$$
\nd This ends the proof of lemma \ref{J_finito}. $\hfill{\rule{2mm}{2mm}}$
\vskip.1cm

\nd \proof {\bf of Lemma \ref {conv_compactos}}  Since $\widetilde{J}$ is finite pick  $\delta>0$ such that $B_{\delta}(x_j)\cap B_{\delta}(x_j)=\emptyset$ for $i\neq j$ with $i, j \in \widetilde{J}$. Next take a compact set $\displaystyle K_{\delta}\subset \r^N \backslash \cup_{j\in\widetilde{J}} B_\delta(x_j)$ and
$\chi\in {C}_0^\infty$ such that
$$
0\leq\chi\leq 1,~~\chi=1~\mbox{on}~K_{\delta},~~\mbox{supp}(\chi)\cap\left(\cup_{j\in\widetilde{J}} B_{\frac{\delta}{2}}(x_j)\right)=\emptyset.
$$
\nd Notice that
$$
\Phi_*(u_n-u)\rightharpoonup\nu~~\mbox{and}~~ \nu = \Phi_*(0)+\sum_{j\in\widetilde{J}}\nu_j\delta_{x_j}~\mbox{in}~ \mathcal{M}.
$$
\nd On the other hand,
$$
 \displaystyle 0\leq\int_{K_{\delta}} \Phi_*(u_n-u)dx  \leq  \displaystyle\int_{\r^N}\Phi_*(u_n-u)\chi dx,
$$
$$
 \displaystyle \int_{\r^N}\Phi_*(u_n-u)\chi dx \rightarrow \int_{\r^N}\chi d\nu,
$$
$$
\int_{\r^N}\chi d\nu  =
 \displaystyle\sum_{j\in\widetilde{J}} \chi(x_j) = 0.
$$
\nd Thus
$$
\int_{K_{\delta}}\Phi_*(u_n-u)dx \rightarrow 0.
$$
\nd Since the argument above holds for each $\delta > 0$ we infer that (\ref{conv-K}) holds for each compact set
$K\subseteq\r^N-\{x_j\}_{j\in\widetilde{J}}$ .
\begin{center}
\section{Proofs of the Main Results}
\end{center}
\begin{lema}\label{convergencia_rho_n_qtp}
    $\rho_n(x)\longrightarrow\rho(x)$~ and~  $\rho(x)\in\partial j(x,u(x))$ a.e.  $x\in\Omega$.
\end{lema}
\nd \proof We will show, at first that
$$
\rho(x)\in\partial j(x,u(x))~ \mbox{a.e.}~ x\in\Omega.
$$
\nd Indeed, let $K\subseteq \r^N \backslash \{x_j\}_{j\in\widetilde{J}}$ be a compact set and take $\varphi\in L_{\widetilde{\Phi}_*}(K)$.
\vskip.1cm

\nd Since
$$
\rho_n\in\partial J(u_n),~~ \rho_n\rightharpoonup \rho~~\mbox{in}~~L_{\widetilde{\Phi}_*}(\Omega)~~ \mbox{and}~~ \rho_n(x)=0~ \mbox{for}~x \in\r^N \backslash \Omega
$$
\nd then
$$
\rho_n\rightharpoonup \rho~~\mbox{in}~~L_{\widetilde{\Phi}_*}(K)
$$
\nd and so
$$
\rho_n\stackrel{*}\rightharpoonup \rho,~~\mbox{em}~~L_{\widetilde{\Phi}_*}(K).
$$
\nd On the other hand, by lemma \ref{conv_compactos},
$$
u_n\stackrel{L_{\Phi_*}(K)}\longrightarrow u
$$
\nd and by \cite[Proposition 2.1.5]{Clarke1}, $\rho\in\partial J(u)$. By the Aubin-Clarke theorem (cf. lemma \ref{aubin_clarke_Le} above ),
$$\rho\in L_{\widetilde{\Phi}_*}(K)~\mbox{and}~\rho(x) \in \partial j(x,u(x))~\mbox{a.e.}~x\in K.
$$
\nd Since
\begin{equation}\label{partition}
\r^N-\{x_j\}_{j\in\widetilde{J}}=\bigcup_{\nu=1}^\infty K_\nu,
\end{equation}
\nd where $\{K_{\nu}\}_{\nu = 1}^{\infty} $ is a sequence of compact sets, it follows that $\rho(x) \in\partial j(x,u(x))$ a.e. $x\in \Omega$.
\vskip.1cm

\nd Next we will show that
$$
\rho_n(x)\longrightarrow\rho(x)~\mbox{a.e.}~x\in \Omega.
$$
\nd Indeed, take  $\varphi_{\nu}\in {C}_0^\infty(\r^N)$ such that $supp(\varphi_{\nu})= K_\nu$. Then
$$
\int_\Omega(\rho_n-\rho)\varphi_\nu dx=\int_{K_\nu}(\rho_n-\rho)\varphi_\nu dx\stackrel{n\rightarrow\infty}\longrightarrow 0,
$$
\nd As a consequence,
$$
(\rho_n-\rho)\varphi_\nu \stackrel{n\rightarrow\infty}\longrightarrow 0~ \mbox{a.e. in}~K_\nu,
$$
\nd so that
$$
\rho_n-\rho\stackrel{n\rightarrow\infty}\longrightarrow 0~ \mbox{a.e. in}~K_\nu.
$$
\nd Therefore
$$
\rho_n-\rho\stackrel{n\rightarrow\infty}\longrightarrow 0~ \mbox{a.e. in }~\r^N.
$$
\nd Since $\rho_n=0$ on $\r^N-\Omega$, it follows that
$$
\rho_n-\rho\stackrel{n\rightarrow\infty}\longrightarrow 0~ \mbox{a.e. in}~\Omega.
$$
\nd This ends the proof of lemma \ref{convergencia_rho_n_qtp}.   $\hfill{\rule{2mm}{2mm}}$
\vskip.1cm

\nd The proof of the next lemma is based on lemma 4.5 in \cite{Fuk_1}.

\begin{lema}\label{conv_grad_qtp}

 ~~~~~~~~   $\nabla u_n(x)~\stackrel{n}\rightarrow~ \nabla u(x)~\mbox{a.e.}~x\in\Omega$.
\end{lema}

\nd \proof Let $\{K_\nu\}_{\nu =1}^{\infty}$ be a family of compact sets such that $(\ref{partition})$ holds. Pick an integer $\nu \geq 1$ and a function $\chi\in {C}_0^\infty(\r^N)$ such that $0\leq\chi\leq 1$,~ $\chi=1~\mbox{on}~K_\nu$ and supp($\chi) \cap\{x_j\}_{j\in\widetilde{J}}\neq\emptyset$.
\vskip.1cm

\nd Set $v_n =\chi (u_n-u)$. It follows that  $v_n$ is bounded in  $\w(\r^N)$ and since  $\langle m(u_n),v_n\rangle) \to 0$
we infer that
\begin{equation}\label{sol_eq}
 \int_{\r^N}\phi(|\nabla u_n|)\nabla u_n\nabla v_ndx-\lambda\int_{\r^N}hv_ndx-\int_{\r^N}\rho_nv_ndx = o_n(1).
\end{equation}
\nd Setting $S_{\chi} = supp(\chi)$ we get
\begin{equation}\label{conv_gradiente-1}
\begin{array}{ll}
  \displaystyle   \displaystyle \int_{S_{\chi}}\phi(|\nabla u_n|)\nabla u_n(\nabla u_n-\nabla u)dx
   + \displaystyle\int_{S_{\chi}}\phi(|\nabla u_n|)\nabla u_n\nabla \chi(u_n- u)dx \\ \\
    =  \displaystyle \int_{S_{\chi}}hv_ndx+\int_{S_{\chi}}\rho_nv_ndx+o_n(1).
\end{array}
\end{equation}
\nd Notice that
$$
\displaystyle \int_{S_{\chi}} \big|\phi(|\nabla u_n|)\nabla u_n\nabla \chi(u_n- u) \big| dx \leq \displaystyle\|\phi(|\nabla u_n|)|\nabla u_n|\|_{L_{\widetilde{\Phi}}(S_{\chi})} |\nabla \chi|_\infty\|(u_n-u)\|_{L_{\Phi}(S_{\chi})} = o_n(1),
$$
$$
\int_{S_{\chi}}hv_ndx= o_n(1),
$$
\nd and since $(\rho_n)$ is bounded in  $L_{\widetilde{\Phi}_*}(\Omega)$,
\[
  \displaystyle \int_{\r^N} \big | \rho_nv_n \big |dx  \leq  \|\rho_n\|_{\widetilde{\Phi}_*}|\chi|_\infty\|(u_n-u)\|_{L_{\widetilde{\Phi}_*}(S_{\chi})}
   = o_{n}(1),
\]
\nd which shows via (\ref{conv_gradiente-1})   that
$$
\int_{K_\nu}\phi(|\nabla u_n|)\nabla u_n(\nabla u_n-\nabla u)dx\stackrel{n} \rightarrow 0.
$$
\nd Using the well known fact that $-\Delta_{\Phi}$ is a map of type $(S_+)$,
$$
\|\nabla u_n-\nabla u\|_{L_{\Phi}(K_\nu)} \stackrel{n} \rightarrow 0.
$$
\nd It follows that
$$
\nabla u_n \stackrel{n}\rightarrow \nabla u ~\mbox{a.e. on}~ K_\nu
$$
\nd and as a consequence,
$$
\nabla u_n \rightarrow \nabla u~\mbox{a.e. on}~ \r^N.
$$
\nd Recalling that $u_n(x)=0~\mbox{for}~x\in \r^N \backslash \Omega$,  we get to
$$
\nabla u_n \rightarrow \nabla u~\mbox{a.e. in}~ \Omega,
$$
\nd endding the proof of  lemma  \ref{conv_grad_qtp}.  $\hfill{\rule{2mm}{2mm}}$
\vskip.2cm
\begin{lema}\label{c_Phi_til}
~~~~~~  $\phi(|\nabla u_n|)\nabla u_n\rightharpoonup\phi(|\nabla u|)\nabla u~~  \mbox{in}~~ \displaystyle \prod  L_{\widetilde{\Phi}}(\Omega).$
\end{lema}

\nd \proof By lemma \ref{conv_grad_qtp},
$$
\nabla u_n  \rightarrow \nabla u~\mbox{a.e. in}~ \Omega.
$$
\nd Since $\phi$ is continuous,
$$
\phi(|\nabla u_n|)\nabla u_n\longrightarrow\phi(|\nabla u|)\nabla u~ \mbox{a.e. in}~\Omega.
$$
\nd Applying lemma 2 in  Gossez \cite[p 88]{gossez-Czech},  ends the proof of lemma \ref{c_Phi_til}.  $\hfill{\rule{2mm}{2mm}}$
\vskip.3cm

\nd \proof {\bf of Theorem \ref{teor_principal}}  By lemma \ref{c_Phi_til},
$$
\int_\Omega\phi(|\nabla u_n|)\nabla u_n\nabla v dx\longrightarrow\int_\Omega\phi(|\nabla u|)\nabla u\nabla v dx,~~v\in\w.
$$
\nd On the other hand,
$$
\int_\Omega\rho_nv dx \longrightarrow \int_\Omega\rho v dx,~~v\in\w,
$$
\nd where
$$
\rho \in L_{\widetilde{\Phi}_*}(\Omega)~~\mbox{and}~~\rho_n(x)\in\partial j(x,u_n(x))~\mbox{a.e.}~x\in\Omega.
$$
\nd Passing to the limt in $(\ref{e1})$ we get to
$$
\int_\Omega\phi(|\nabla u |)\nabla u \nabla v dx-\lambda\int_\Omega hv dx-\int_\Omega \rho vdx=0,~~v\in\w.
$$
\nd Thus  $u \in\w$ is a solution of (\ref{problema}), in the sense of  Definition \ref{solutionDef-a} and since
$h \ne 0$,  we get   $u \not\equiv 0$.
\vskip.2cm

\nd \nd {\bf Claim.} $u \geq0$. Indeed, note that
$$
u_n=u_n^+ - u_n^-,~\nabla u_n=\nabla u_n^+ - \nabla u_n^-~\mbox{and}~|\nabla u_n|^2=|\nabla u_n^+|^2+|\nabla u_n^-|^2
$$
\nd Thus
$$
\begin{array}{lll}
  \displaystyle\int_\Omega\Phi(|\nabla u_n^-|)dx & \leq & \displaystyle\int_\Omega\Phi([|\nabla u_n^-|^2+|\nabla u_n^+|^2]^{\frac{1}{2}})dx \\ \\
   & = & \displaystyle \int_\Omega\Phi(|\nabla u_n|)dx
\end{array}
$$
\nd so that $(u_n^-)$ is bounded in $\w$. Noting that $\langle w_n, -u_n^-\rangle = o_n(1)$ we have
$$
\begin{array}{lcl}
 o_n(1) &=& \displaystyle - \int_\Omega\phi(|\nabla u_n|)\nabla u_n\nabla u_n^-dx + \lambda\int_\Omega hu_n^-dx  +\int_\Omega\rho_nu_n^-dx \\
    &=& \displaystyle \int_\Omega\phi(|\nabla u_n^-|)|\nabla u_n^-|^2dx  + \lambda\int_\Omega hu_n^-dx  + \int_\Omega\rho_nu_n^-dx   \\
    &\geq&  \ell \displaystyle \int_\Omega\Phi(|\nabla u_n^-|)dx.
\end{array}
$$
\nd Thus
$$
 \displaystyle \int_\Omega\Phi(|\nabla u_n^-|)dx \rightarrow 0,
$$
\nd and hence  $u_n^-\rightarrow 0~~\mbox{in}~~\w$, showing that $u \geq 0$.
\vskip.2cm

\nd \proof  {\bf of (\ref{maineq ae})}~ Since $u$ is a solution of (\ref{problema}), there is $\rho := \rho_u \in  L_{{\widetilde{\Phi}_*}}(\Omega)$ such that
 \[
        \int_\Omega\phi(|\nabla u|)\nabla u\nabla vdx=\int_\Omega\rho vdx+\lambda\int_\Omega hvdx,~v\in C_0^{\infty}(\Omega).
   \]
\nd By Remark \ref{orlicz-sobolev}, $\Delta_\Phi u \in L_{\widetilde{\Phi}_*}(\Omega~)$. Since also $h \in L_{\widetilde{\Phi}_*}(\Omega~)$ it follows that
 \[
        \int_\Omega [-\Delta_\Phi u~ - \rho - \lambda h]vdx = 0,~v\in C_0^{\infty}(\Omega).
   \]
\nd Hence
$$
-\Delta_{\Phi} u = \rho + \lambda h~~ \mbox{a.e. in}~~ \Omega.
$$
$ \hfill{\rule{2mm}{2mm}}$

\begin{center}
\section{ Appendix }
\end{center}

\nd The results below are elementary and can be found in    \cite{Fuk_1, Fuk_2}.

\begin{lema}\label{lema_naru}
        Assume  $(\phi_1)-(\phi_3)$. Let
        $$
\zeta_0(t)=\min\{t^\ell,t^m\}~~\mbox{and}~~\zeta_1(t)=\max\{t^\ell,t^m\},~~t\geq 0.
$$
  \nd Then
        $$
            \zeta_0(\rho)\Phi(t)\leq\Phi(\rho t)\leq \zeta_1(\rho)\Phi(t),~\rho, t> 0,
        $$
        $$
            \zeta_0(\|u\|_{\Phi})\leq\int_\Omega\Phi(u)dx\leq \zeta_1(\|u\|_{\Phi}),~ u\in L_{\Phi}(\Omega).
        $$
\end{lema}

\begin{lema}\label{lema_naru_1}
        Assume  $(\phi_1)-(\phi_3)$. Let
        $$\zeta_2(t)=\min\{t^{\ell^*},t^{m^*}\}~~\mbox{and}~~\zeta_2(t)=\max\{t^{\ell^*},t^{m^*}\},~~t\geq 0.$$
  \nd Then
        $$
            \zeta_2(\rho)\Phi_*(t)\leq\Phi_*(\rho t)\leq \zeta_3(\rho)\Phi_*(t),~\rho, t> 0,
        $$
        $$
            \zeta_2(\|u\|_{\Phi_*})\leq\int_\Omega\Phi_*(u)dx\leq \zeta_3(\|u\|_{\Phi_*}),~ u\in L_{\Phi_*}(\Omega).
        $$
\end{lema}

\begin{center}

\end{center}

\begin{flushright}
{J. V. Gon\c{c}alves}\\
\smallskip
  \scriptsize{Universidade Federal de Goi\'as\\
   Instituto de Matem\'atica e Estat\'istica\\
   74001-970 Goi\^ania, GO - Brasil}
\end{flushright}

\begin{flushright}
{M. L. Carvalho}\\
\scriptsize{Universidade Federal de Goi\'as\\
 Departamento de Matem\'atica\\
 75804-020  Jata\'i, GO - Brasil}
\end{flushright}

\end{document}